\documentclass [12pt]{article}
\usepackage{amsmath}
\usepackage{amsfonts}
\usepackage{amssymb,latexsym}
\newcommand{\qed}{\ifmmode$\Box$\else{\unskip\nobreak\hfil
		\penalty50\hskip1em\null\nobreak\hfil$\Box$
		\parfillskip=0pt\finalhyphendemerits=0\endgraf}\fi}
\begin{document}
\title
{
On the Moisil-Theodoresco operator in orthogonal curvilinear coordinates.
}
\title
{
	On the Moisil-Theodoresco operator in orthogonal curvilinear coordinates
}
\author{Juan Bory-Reyes$^{(1)}$ and Marco Antonio P\'{e}rez-de la Rosa$^{(2)}$}

\date{\small $^{(1)}$ ESIME-Zacatenco. Instituto Polit\'ecnico Nacional. CD-MX. 07738. M\'exico.\\E-mail: juanboryreyes@yahoo.com\\
	$^{(2)}$ Department of Actuarial Sciences, Physics and Mathematics, Universidad de las Am\'{e}ricas Puebla.
	San Andr\'{e}s Cholula, Puebla. 72810. M\'{e}xico.\\ Email: marco.perez@udlap.mx}
\maketitle

\begin{abstract}
	It is generally well understood the legitimate action of the Moisil-Theo\-do\-res\-co ope\-ra\-tor, over a quaternionic valued function defined on $\mathbb{R}^3$ (sum of a scalar and a vector field) in Cartesian coordinates, but it does not so in any orthogonal curvilinear coordinate system. This paper sheds some new light on the technical aspect of the subject. Moreover, we introduce a notion of quaternionic Laplace operator acting on a quaternionic valued function from which one can recover both scalar and vector Laplacians in vector analysis context.
\end{abstract}
{\bf Keywords.} Moisil-Theodoresco operator \and Laplace operator \and hyperholomorphic functions \and orthogonal curvilinear coordinates\\
{\bf MSC (2000).} Primary 30G35 \and Secondary 35J05

\section{Introduction.}

The Moisil-Theodoresco operator (a determined first-order elliptic system to be defined below) is nowadays considered to be a good analogue of the usual Cauchy-Riemann operator of complex analysis to the quaternionic setting. It is a square root of the scalar Laplace operator in $\mathbb{R}^3$. The consideration of the Moisil-Theodoresco operator was the starting point in the development of the hyperholomorphic function theory, see \cite{Fu} and \cite{MT}.

In recent times, the theory of hyperholomorphic quaternion-valued functions has been developed along many interesting directions, including
the treatment of several problems from mathematical physics with quaternionic analysis techniques, in particular the Lam\'{e}-Navier system from the theory of elasticity. See for instance \cite{Grigo1}, \cite{Grigo2}, \cite{Grigo3}, \cite{GuSp} and \cite{MoMoAbBo}.

It is generally well understood the legitimate action of the Moisil-Theo\-do\-res\-co operators, over a quaternionic valued function defined on $\mathbb{R}^3$ (sum of a scalar and a vector field) in Cartesian coordinates, but it does not so in any orthogonal curvilinear coordinate system. This work is intended to derive the explicit expressions for such operator in a general orthogonal curvilinear coordinate system.

It is worth pointing out that the product of the Moisil-Theodoresco operator with its adjoint produces and elliptic (but not strongly elliptic) operator in $\mathbb{R}^3$, which can be see as a generalization of the Bitsadze operator to the space $\mathbb{R}^3$, see \cite{Bit}, \cite{Dz} and \cite{S}.

The term quaternionic analysis is a generalization or extension of complex analysis for holomorphic functions depending on one complex variable, where the latter means function theory for complex-valued functions in one complex variable.

The class of functions belonging to the kernel of the Moisil-Theodoresco operator enjoys many satisfactory properties and strategies similarly to holomorphic functions with results known from vector analysis.

Vector analysis is a classical subject dealing with those aspects of vectors which are usually confined to three dimensional Euclidean space.  Vector multiplication is not uniquely defined, such as the scalar and vector products. The last is not associative without unit.

It is customary to think of a pure quaternion as a vector, although there are reasons why this is quite wrong. In quaternionic system the multiplication conjugates both the scalar and vector products, meanwhile in vector system they remain independently. For a deeper discussion of the differences in the multiplication proprieties of two pure quaternions and two vectors we refer the reader to \cite{Ma}.

The study of quaternionic analysis in a Euclidean setting and Cartesian coordinates has grown into several major research fields, including many applications to other branches of mathematics, physical sciences, and engineering. For a thorough treatment of the subject can be found in \cite{Krav} and \cite{KravShap} and the references given there.

At the same time, a closer examination reveals that, surprisingly for us, the definite meaning of the Moisil-Theodoresco operator acting on a quaternionic valued function defined on $\mathbb{R}^3$ (sum of a scalar and a vector fields) in any orthogonal curvilinear coordinate system greatly differ from an approach within the framework of Cartesian coordinates. The purpose of the present work is to describe these subtle differences. The authors are not aware that the results of this paper can be found in the literature.

For the applications, we consider the relation of the Laplace operator and the product of the Moisil-Theodoresco operator with itself and so with its adjoint to the equilibrium of an isotropic elastic body in $\mathbb{R}^3$. Along the way a formula for the Lam\'{e} operator in any orthogonal curvilinear coordinate system is obtained.

\section{Preliminaries}
In this section, we provide standard facts from classical quaternionic analysis to be used in this paper. For more information, we refer the reader to \cite{KravShap}.

Let $\mathbb H(\mathbb C)$ be the set of complex quaternions, i.e., that each quaternion $a$ is represented in the form $a = \sum_{k=0}^3 a_k\mathbf{i}_k$, with $\{a_k\}\subset\mathbb C$; $\mathbf{i}_0 = 1$ stands for the unit and $\mathbf{i}_1, \mathbf{i}_2, \mathbf{i}_3$ stand for the quaternionic imaginary units.
Denote the complex imaginary unit in $\mathbb{C}$ by $i$ as usual. By definition, $i$ commutes with all the quaternionic imaginary units $\mathbf{i}_1, \mathbf{i}_2, \mathbf{i}_3$.

The set $\mathbb H(\mathbb C)$ is a complex non-commutative, associative algebra with zero divisors. The involution $a\rightarrow \bar{a}$, called quaternionic conjugation, is defined by
\[
\bar{a} := \sum^3_{k=0} a_k\cdot\bar{\mathbf{i}}_k = a_0 - \sum^3_{k=1} a_k\cdot \mathbf{i}_k.
\]
It satisfies $\overline{ab} = \bar{b}\bar{a}$. The Euclidean norm $\lvert a\rvert$ in $\mathbb H(\mathbb C)$ is defined by $\lvert a\rvert := \sqrt{\sum_{k=0}^3\lvert a_k\rvert^2}$. 

Writing for $a = \sum_{k=0}^3 a_k\mathbf{i}_k\in\mathbb H(\mathbb C)$, $a_0 =: \text{Sc}(a)$, $\vec{a} := \sum_{k=1}^3 a_k\mathbf{i}_k =: \text{Vec}(a)$, we have $a = a_0 + \vec{a}$. We call $a_0$ the scalar part of the complex quaternion $a$ and $\vec{a}$ the vector part of $a$. Then $\{ \text{Vec}(a) : a\in\mathbb H(\mathbb C)\}$ is identified with $\mathbb C^3$. This enables us to write $\bar a = \text{Sc}(a) - \vec a$.

For any $a,b\in\mathbb H(\mathbb C)$:
\[
a\,b := a_0\,b_0 - \langle\vec a,\vec b\rangle + a_0\,\vec b + b_0\,\vec a + [\vec a,\vec b],
\]
where
\[
\langle\vec a,\vec b\rangle := \sum_{k=1}^3 a_k\,b_k,\,\, [\vec a,\vec b]:=
\left |
\begin{array}{rrr}
\mathbf{i}_1 & \mathbf{i}_2 & \mathbf{i}_3\\
a_1& a_2 & a_3\\
b_1& b_2 & b_3
\end{array}
\right |.
\]
In particular, if $a_0=b_0=0$ then $a\,b :=  - \langle\vec a,\vec b\rangle + [\vec a,\vec b]$.

Denote by $\mathfrak{S}$ the set of zero divisors from $\mathbb{H}(\mathbb{C})$ and by $G\mathbb{H}(\mathbb{C})$ the subset of in\-ver\-ti\-ble elements from $\mathbb{H}(\mathbb{C})$. If $a\notin\mathfrak{S}\cup\{0\}$ then $\displaystyle a^{-1}:=\frac{\bar{a}}{(a\bar{a})} $ is the inverse of the complex quaternion $a$. Note that $G\mathbb{H}(\mathbb{C})=\mathbb{H}(\mathbb{C})\setminus(\mathfrak{S}\cup\{0\})$.

Let $\Omega$ be a domain in $\mathbb{R}^3$. On $C^1(\Omega;\mathbb{H}(\mathbb{C}))$ the Moisil-Theodoresco operator, denoted by $D_{MT}$, is defined to be:
\begin{equation}
D_{MT}[f]:=\mathbf{i}_1\frac{\partial f}{\partial x}+\mathbf{i}_2\frac{\partial f}{\partial y}+\mathbf{i}_3\frac{\partial f}{\partial z}.
\end{equation}

Notice that $D_{MT}$ factorizes the scalar Laplacian as follows:
\begin{equation}
-D_{MT}^2=\Delta_{\mathbb{R}^4},
\end{equation}
which implies several advantages in the applications to physical problems.

The operator $\Delta_{\mathbb{R}^4}$ is a scalar operator, its acts separately on every coordinate function $f_k$ of $f$ as
$$\Delta_{\mathbb{R}^4}[f]:=\Delta[f_0]+\mathbf{i}_1\Delta[f_1]+\mathbf{i}_2\Delta[f_2]+\mathbf{i}_3\Delta[f_3].$$

This property guarantees that any hyperholomorphic function is also harmonic.

Note that for an $\mathbb{H}(\mathbb{C})$-valued function $f:=f_0+ \vec{f}$ the action of the operator $D_{MT}$ can be
represented as follows
\begin{equation}
D_{MT}[f]=-\mathrm{div}[\vec{f}]+\mathrm{grad}[f_0]+\mathrm{curl}[\vec{f}].
\end{equation}
This is an immediate consequence of the quaternionic product.

As usual, the Moisil-Theodoresco operator can act on the right in which case the notation $D_{MT}^r[f]$  for the same $f=f_0+\vec f$ means that:
\begin{equation}
D_{MT}^r[f]:=\frac{\partial f}{\partial x}\mathbf{i}_1+\frac{\partial f}{\partial y}\mathbf{i}_2+\frac{\partial f}{\partial z}\mathbf{i}_3.
\end{equation}

Furthermore, consider the determined non-homogeneous first-order elliptic system
\begin{equation}\label{nhMT}
\begin{cases}
\mathrm{div}[\vec{f}]=g_0,\\
\mathrm{grad}[f_0]+\mathrm{curl}[\vec{f}]=\vec{g},
\end{cases}
\end{equation}
and also its adjoint
\begin{equation}\label{adjMT}
\begin{cases}
\mathrm{div}[\vec{f}]=0,\\
\mathrm{grad}[f_0]-\mathrm{curl}[\vec{f}]=0,
\end{cases}
\end{equation}
which correspond to the Moisil-Theodoresco operator $D_{MT}$ and $D_{MT}^r$, respectively. If $(f_0, \vec{f})$ associates with the $\mathbb{H}(\mathbb{C})$-valued function $f:=f_0+ \vec{f}$, then $(f_0, \vec{f})$ satisfies systems (\ref{nhMT}) and (\ref{adjMT}) are equivalent to $f$ satisfying the equations
$$D_{MT}[f]=g,$$
where $g:=g_0+ \vec{g}$ and
$$D_{MT}^r[f]=0,$$
respectively. A. Dzhuraev considers in his book \cite{Dz} the following matrix differential operator named also Moisil-Theodoresco operator
\begin{equation}\label{repr_matr_MT}
\left(\begin{array}{cccc}
    0 & \mathrm{div} \\
    \mathrm{grad} & \mathrm{curl}
  \end{array}\right),
\end{equation}
whose action can be identify with the negative quaternionic conjugated to ${D}_{MT}.$
\section{The Moisil-Theodoresco operator in general orthogonal curvilinear coordinates}
Let $x,y$ and $z$ represent the Cartesian coordinates in $\mathbb{R}^3$, with the corresponding orthonormal basis $\{\mathbf{i}_1,\mathbf{i}_2,\mathbf{i}_3\}$. We assume that $q_1, q_2, q_3$ are general orthogonal curvilinear coordinates defined via three invertible and continuously differentiable transformation functions
\[
x=x(q_1,q_2,q_3),\;\; y=y(q_1,q_2,q_3),\;\; z=z(q_1,q_2,q_3).
\]
Each system is characterized by the metric coefficients:
\begin{equation}
h_i=\sqrt{\left(\frac{\partial x}{\partial q_i}\right)^2+\left(\frac{\partial y}{\partial q_i}\right)^2+\left(\frac{\partial z}{\partial q_i}\right)^2},\;\;\;i=1,2,3,
\end{equation}
and the unit vectors $\mathbf{u}_1$ , $\mathbf{u}_2$ , and $\mathbf{u}_3$ can be obtained as follows:

\begin{equation}
\mathbf{u}_i=\frac{1}{h_i}\left(\frac{\partial x}{\partial q_i}\,\mathbf{i}_1+\frac{\partial y}{\partial q_i}\,\mathbf{i}_2+\frac{\partial z}{\partial q_i}\,\mathbf{i}_3\right),\;\;\;i=1,2,3.
\end{equation}

By abuse of notation, we continue to write a quaternionic-valued function $f$ in orthogonal curvilinear coordinates in the following way:
\begin{equation}
f=f_0+\vec{f}=f_0+f_1\,\mathbf{u}_1+f_2\,\mathbf{u}_2+f_3\,\mathbf{u}_3.
\end{equation}

It is known (see, e.g. \cite{Bar}, \cite{Gri}, \cite{Mo}, \cite{Mo2} and \cite{Tai}) that the vectorial operations gradient, divergence and curl in orthogonal curvilinear coordinates can be expressed in terms of the metric coefficients and take the following form:
\begin{equation}
\mathrm{grad}[f_0]=\frac{1}{h_1}\frac{\partial f_0}{\partial q_1}\,\mathbf{u}_1+\frac{1}{h_2}\frac{\partial f_0}{\partial q_2}\,\mathbf{u}_2+\frac{1}{h_3}\frac{\partial f_0}{\partial q_3}\,\mathbf{u}_3,
\end{equation}

\begin{equation}
\mathrm{div}[\vec{f}]=\frac{1}{h_1h_2h_3}\left[\frac{\partial (h_2h_3f_1)}{\partial q_1}+\frac{\partial (h_1h_3f_2)}{\partial q_2}+\frac{\partial (h_1h_2f_3)}{\partial q_3}\right],
\end{equation}

and

\begin{align}
\mathrm{curl}[\vec{f}]=&\frac{1}{h_2h_3}\left[\frac{\partial (h_3f_3)}{\partial q_2}-\frac{\partial (h_2f_2)}{\partial q_3}\right]\,\mathbf{u}_1+\frac{1}{h_1h_3}\left[\frac{\partial (h_1f_1)}{\partial q_3}-\frac{\partial (h_3f_3)}{\partial q_1}\right]\,\mathbf{u}_2\notag\\
&+\frac{1}{h_1h_2}\left[\frac{\partial (h_2f_2)}{\partial q_1}-\frac{\partial (h_1f_1)}{\partial q_2}\right]\,\mathbf{u}_3.
\end{align}
Then, the Moisil-Theodoresco operator for a quaternionic-valued function in orthogonal curvilinear coordinates is defined by:
\begin{align}
D_{MT}[f]&=-\mathrm{div}[\vec{f}]+\mathrm{grad}[f_0]+\mathrm{curl}[\vec{f}]\notag\\
&=-\frac{1}{h_1h_2h_3}\left[\frac{\partial (h_2h_3f_1)}{\partial q_1}+\frac{\partial (h_1h_3f_2)}{\partial q_2}+\frac{\partial (h_1h_2f_3)}{\partial q_3}\right]\notag\\
&\;\;\;\;+\left\{\frac{1}{h_1}\frac{\partial f_0}{\partial q_1}+\frac{1}{h_2h_3}\left[\frac{\partial (h_3f_3)}{\partial q_2}-\frac{\partial (h_2f_2)}{\partial q_3}\right]\right\}\,\mathbf{u}_1,\notag\\
&\;\;\;\;+\left\{\frac{1}{h_2}\frac{\partial f_0}{\partial q_2}+\frac{1}{h_1h_3}\left[\frac{\partial (h_1f_1)}{\partial q_3}-\frac{\partial (h_3f_3)}{\partial q_1}\right]\right\}\,\mathbf{u}_2,\notag\\
&\;\;\;\;+\left\{\frac{1}{h_3}\frac{\partial f_0}{\partial q_3}+\frac{1}{h_1h_2}\left[\frac{\partial (h_2f_2)}{\partial q_1}-\frac{\partial (h_1f_1)}{\partial q_2}\right]\right\}\,\mathbf{u}_3.
\end{align}

\subsection{Scalar versus vector Laplacian}
The action of the vector Laplacian is reminiscent of the action of scalar Laplacian. Whereas the scalar Laplacian applies to a scalar field and returns a scalar quantity, the vector Laplacian applies to a vector field, returning a vector quantity. In orthonormal Cartesian coordinates, the returned vector field is equal to the vector field of the scalar Laplacian acting on each Cartesian component of the vector field separately. Meanwhile, in general orthogonal curvilinear coordinates they are two entirely different operators. In recent developments, authors (see \cite{HiChi}, \cite{Ken}, \cite{Mo}, \cite{Mo2}, \cite{Red} and \cite{Tai}) have been studied extensively the general expressions of the vector Laplacian in in orthogonal curvilinear coordinates.

In this way vector analysis in general curvilinear coordinates systems contains two Laplacian operators neither of which can be applied to either scalar or vector fields. At the same time quaternion-valued functions defined in $\mathbb{R}^3$ can be identified with a sum of a scalar and vector field, what would be impossible in vector calculus. Based on this fact naturally raise the question of whether or not a quaternionic Laplacian could be defined. In that respect we define on $C^2(\Omega;\mathbb{H}(\mathbb{C}))$ the following operator
\begin{equation}
\Delta_{\mathbb H}[f]:=\Delta_0[f_0]+\vec{\Delta}[\vec{f}],
\end{equation}
where
\begin{align}
\Delta_0[f_0]:&=\mathrm{div}[\mathrm{grad}[f_0]]\notag\\
&=\frac{1}{h_1h_2h_3}\left[\frac{\partial}{\partial q_1}\left(\frac{h_2h_3}{h_1}\frac{\partial f_0}{\partial q_1}\right)+\frac{\partial}{\partial q_2}\left(\frac{h_1h_3}{h_2}\frac{\partial f_0}{\partial q_2}\right)+\frac{\partial}{\partial q_3}\left(\frac{h_1h_2}{h_3}\frac{\partial f_0}{\partial q_3}\right)\right],
\end{align}
and
\begin{align}
\vec{\Delta}[\vec{f}]:&=\mathrm{grad}[\mathrm{div}[\vec{f}]]-\mathrm{curl}[\mathrm{curl}[\vec{f}]]\notag\\
&=\left\{\frac{1}{h_1}\frac{\partial}{\partial q_1}\left[\frac{1}{h_1h_2h_3}\frac{\partial(h_2h_3f_1)}{\partial q_1}\right]+\frac{1}{h_1}\frac{\partial}{\partial q_1}\left[\frac{1}{h_1h_2h_3}\frac{\partial(h_1h_3f_2)}{\partial q_2}\right]\right.\notag\\
&\;\;\;\;\;+\frac{1}{h_1}\frac{\partial}{\partial q_1}\left[\frac{1}{h_1h_2h_3}\frac{\partial(h_1h_2f_3)}{\partial q_3}\right]-\frac{1}{h_2h_3}\frac{\partial}{\partial q_2}\left[\frac{h_3}{h_1h_2}\frac{\partial(h_2f_2)}{\partial q_1}\right]\notag\\
&\;\;\;\;\;+\frac{1}{h_2h_3}\frac{\partial}{\partial q_2}\left[\frac{h_3}{h_1h_2}\frac{\partial(h_1f_1)}{\partial q_2}\right]+\frac{1}{h_2h_3}\frac{\partial}{\partial q_3}\left[\frac{h_2}{h_1h_3}\frac{\partial(h_1f_1)}{\partial q_3}\right]\notag\\
&\;\;\;\;\;\left.-\frac{1}{h_2h_3}\frac{\partial}{\partial q_3}\left[\frac{h_2}{h_1h_3}\frac{\partial(h_3f_3)}{\partial q_1}\right]\right\}\,\mathbf{u}_1\notag\\
\;\;\;\;\;\;\;\;\;\;&\;\;\;+\left\{\frac{1}{h_2}\frac{\partial}{\partial q_2}\left[\frac{1}{h_1h_2h_3}\frac{\partial(h_2h_3f_1)}{\partial q_1}\right]+\frac{1}{h_2}\frac{\partial}{\partial q_2}\left[\frac{1}{h_1h_2h_3}\frac{\partial(h_1h_3f_2)}{\partial q_2}\right]\right.\notag\\
&\;\;\;\;\;+\frac{1}{h_2}\frac{\partial}{\partial q_2}\left[\frac{1}{h_1h_2h_3}\frac{\partial(h_1h_2f_3)}{\partial q_3}\right]-\frac{1}{h_1h_3}\frac{\partial}{\partial q_3}\left[\frac{h_1}{h_2h_3}\frac{\partial(h_3f_3)}{\partial q_2}\right]\notag\\
&\;\;\;\;\;+\frac{1}{h_1h_3}\frac{\partial}{\partial q_3}\left[\frac{h_1}{h_2h_3}\frac{\partial(h_2f_2)}{\partial q_3}\right]+\frac{1}{h_1h_3}\frac{\partial}{\partial q_1}\left[\frac{h_3}{h_1h_2}\frac{\partial(h_2f_2)}{\partial q_1}\right]\notag\\
&\;\;\;\;\;\left.-\frac{1}{h_1h_3}\frac{\partial}{\partial q_1}\left[\frac{h_3}{h_1h_2}\frac{\partial(h_1f_1)}{\partial q_2}\right]\right\}\,\mathbf{u}_2\notag\\
\;\;\;\;\;\;\;\;\;\;&\;\;\;+\left\{\frac{1}{h_3}\frac{\partial}{\partial q_3}\left[\frac{1}{h_1h_2h_3}\frac{\partial(h_2h_3f_1)}{\partial q_1}\right]+\frac{1}{h_3}\frac{\partial}{\partial q_3}\left[\frac{1}{h_1h_2h_3}\frac{\partial(h_1h_3f_2)}{\partial q_2}\right]\right.\notag\\
&\;\;\;\;\;+\frac{1}{h_3}\frac{\partial}{\partial q_3}\left[\frac{1}{h_1h_2h_3}\frac{\partial(h_1h_2f_3)}{\partial q_3}\right]-\frac{1}{h_1h_2}\frac{\partial}{\partial q_1}\left[\frac{h_2}{h_1h_3}\frac{\partial(h_1f_1)}{\partial q_3}\right]\notag\\
&\;\;\;\;\;+\frac{1}{h_1h_2}\frac{\partial}{\partial q_1}\left[\frac{h_2}{h_1h_3}\frac{\partial(h_3f_3)}{\partial q_1}\right]+\frac{1}{h_1h_2}\frac{\partial}{\partial q_2}\left[\frac{h_1}{h_2h_3}\frac{\partial(h_3f_3)}{\partial q_2}\right]\notag\\
&\;\;\;\;\;\left.-\frac{1}{h_1h_2}\frac{\partial}{\partial q_2}\left[\frac{h_1}{h_2h_3}\frac{\partial(h_2f_2)}{\partial q_3}\right]\right\}\,\mathbf{u}_3,\label{VectorDelta}
\end{align}
Notice that the operator $\vec{\Delta}$ is independent of the choice of the coordinate system. In this way we obtain what we shall call the quaternionic Laplace operator in general orthogonal curvilinear coordinates.

A short calculation shows that the Moisil-Theodoresco operator factorize the quaternionic Laplacian in orthogonal curvilinear coordinates as follows:
\begin{equation}
D_{MT}^2[f]=-\Delta_{\mathbb H}[f].
\end{equation}

On the other hand, we introduce an elliptic (but not strongly elliptic) operator in $\mathbb{R}^3$, which can be see as a generalization of the Bitsadze operator to the space $\mathbb{R}^3$ (see \cite{Bit}, \cite{Dz} and \cite{S}) as the vector part of the following opertator:
\begin{equation}
\widetilde{\Delta_{\mathbb H}}[f]:=\Delta_0[f_0]+\widetilde{\vec{\Delta}}[\vec{f}],
\end{equation}
where
\begin{align}
\widetilde{\vec{\Delta}}[\vec{f}]:&=\mathrm{grad}[\mathrm{div}[\vec{f}]]+\mathrm{curl}[\mathrm{curl}[\vec{f}]]\notag\\
&=\left\{\frac{1}{h_1}\frac{\partial}{\partial q_1}\left[\frac{1}{h_1h_2h_3}\frac{\partial(h_2h_3f_1)}{\partial q_1}\right]+\frac{1}{h_1}\frac{\partial}{\partial q_1}\left[\frac{1}{h_1h_2h_3}\frac{\partial(h_1h_3f_2)}{\partial q_2}\right]\right.\notag\\
&\;\;\;\;\;+\frac{1}{h_1}\frac{\partial}{\partial q_1}\left[\frac{1}{h_1h_2h_3}\frac{\partial(h_1h_2f_3)}{\partial q_3}\right]+\frac{1}{h_2h_3}\frac{\partial}{\partial q_2}\left[\frac{h_3}{h_1h_2}\frac{\partial(h_2f_2)}{\partial q_1}\right]\notag\\
&\;\;\;\;\;-\frac{1}{h_2h_3}\frac{\partial}{\partial q_2}\left[\frac{h_3}{h_1h_2}\frac{\partial(h_1f_1)}{\partial q_2}\right]-\frac{1}{h_2h_3}\frac{\partial}{\partial q_3}\left[\frac{h_2}{h_1h_3}\frac{\partial(h_1f_1)}{\partial q_3}\right]\notag\\
&\;\;\;\;\;\left.+\frac{1}{h_2h_3}\frac{\partial}{\partial q_3}\left[\frac{h_2}{h_1h_3}\frac{\partial(h_3f_3)}{\partial q_1}\right]\right\}\,\mathbf{u}_1\notag\\
\;\;\;\;\;\;\;\;\;\;&\;\;\;+\left\{\frac{1}{h_2}\frac{\partial}{\partial q_2}\left[\frac{1}{h_1h_2h_3}\frac{\partial(h_2h_3f_1)}{\partial q_1}\right]+\frac{1}{h_2}\frac{\partial}{\partial q_2}\left[\frac{1}{h_1h_2h_3}\frac{\partial(h_1h_3f_2)}{\partial q_2}\right]\right.\notag\\
&\;\;\;\;\;+\frac{1}{h_2}\frac{\partial}{\partial q_2}\left[\frac{1}{h_1h_2h_3}\frac{\partial(h_1h_2f_3)}{\partial q_3}\right]+\frac{1}{h_1h_3}\frac{\partial}{\partial q_3}\left[\frac{h_1}{h_2h_3}\frac{\partial(h_3f_3)}{\partial q_2}\right]\notag\\
&\;\;\;\;\;-\frac{1}{h_1h_3}\frac{\partial}{\partial q_3}\left[\frac{h_1}{h_2h_3}\frac{\partial(h_2f_2)}{\partial q_3}\right]-\frac{1}{h_1h_3}\frac{\partial}{\partial q_1}\left[\frac{h_3}{h_1h_2}\frac{\partial(h_2f_2)}{\partial q_1}\right]\notag\\
&\;\;\;\;\;\left.+\frac{1}{h_1h_3}\frac{\partial}{\partial q_1}\left[\frac{h_3}{h_1h_2}\frac{\partial(h_1f_1)}{\partial q_2}\right]\right\}\,\mathbf{u}_2\notag\\
\;\;\;\;\;\;\;\;\;\;&\;\;\;+\left\{\frac{1}{h_3}\frac{\partial}{\partial q_3}\left[\frac{1}{h_1h_2h_3}\frac{\partial(h_2h_3f_1)}{\partial q_1}\right]+\frac{1}{h_3}\frac{\partial}{\partial q_3}\left[\frac{1}{h_1h_2h_3}\frac{\partial(h_1h_3f_2)}{\partial q_2}\right]\right.\notag\\
&\;\;\;\;\;+\frac{1}{h_3}\frac{\partial}{\partial q_3}\left[\frac{1}{h_1h_2h_3}\frac{\partial(h_1h_2f_3)}{\partial q_3}\right]+\frac{1}{h_1h_2}\frac{\partial}{\partial q_1}\left[\frac{h_2}{h_1h_3}\frac{\partial(h_1f_1)}{\partial q_3}\right]\notag\\
&\;\;\;\;\;-\frac{1}{h_1h_2}\frac{\partial}{\partial q_1}\left[\frac{h_2}{h_1h_3}\frac{\partial(h_3f_3)}{\partial q_1}\right]-\frac{1}{h_1h_2}\frac{\partial}{\partial q_2}\left[\frac{h_1}{h_2h_3}\frac{\partial(h_3f_3)}{\partial q_2}\right]\notag\\
&\;\;\;\;\;\left.+\frac{1}{h_1h_2}\frac{\partial}{\partial q_2}\left[\frac{h_1}{h_2h_3}\frac{\partial(h_2f_2)}{\partial q_3}\right]\right\}\,\mathbf{u}_3,\label{VectorDeltaTilde}
\end{align}

It turns out that it is also possible to factorize the operator $\widetilde{\Delta_{\mathbb H}}$ in the following manner:
\begin{equation}
D_{MT} D_{MT}^r[f]=-\widetilde{\Delta_{\mathbb H}}[f].
\end{equation}
As we will see below, this operator is related to the equilibrium of an isotropic elastic body (see \cite{Mi}).
\section{Examples}
\subsection{Cartesian coordinates}
For Cartesian coordinates one has $h_1=h_2=h_3=1$ and $\mathbf{u}_1=\mathbf{i}_1,\;\mathbf{u}_2=\mathbf{i}_2\;\mathbf{u}_3=\mathbf{i}_3$.
Hence, the Moisil-Theodoresco operator takes the form:
\begin{align}
D_{MT}[f]&=-\left\{\frac{\partial f_1}{\partial x}+\frac{\partial f_2}{\partial y}+\frac{\partial f_3}{\partial z}\right\}\notag\\
&\;\;\;\;+\left\{\frac{\partial f_0}{\partial x}+\frac{\partial f_3}{\partial y}-\frac{\partial f_2}{\partial z}\right\}\,\mathbf{i}_1\notag\\
&\;\;\;\;+\left\{\frac{\partial f_0}{\partial y}+\frac{\partial f_1}{\partial z}-\frac{\partial f_3}{\partial x}\right\}\,\mathbf{i}_2\notag\\
&\;\;\;\;+\left\{\frac{\partial f_0}{\partial z}+\frac{\partial f_2}{\partial x}-\frac{\partial f_1}{\partial y}\right\}\,\mathbf{i}_3.
\end{align}

The scalar Laplacian is:
\begin{equation}
\Delta_0[f_0]=\frac{\partial^2 f_0}{\partial x^2}+\frac{\partial^2 f_0}{\partial y^2}+\frac{\partial^2 f_0}{\partial z^2},
\end{equation}
and the vector Laplacian:
\begin{align}
\vec{\Delta}[\vec{f}]&=\left\{\frac{\partial^2 f_1}{\partial x^2}+\frac{\partial^2 f_1}{\partial y^2}+\frac{\partial^2 f_1}{\partial z^2}\right\}\,\mathbf{i}_1\notag\\
&\;\;\;+\left\{\frac{\partial^2 f_2}{\partial x^2}+\frac{\partial^2 f_2}{\partial y^2}+\frac{\partial^2 f_2}{\partial z^2}\right\}\,\mathbf{i}_2\notag\\
&\;\;\;+\left\{\frac{\partial^2 f_3}{\partial x^2}+\frac{\partial^2 f_3}{\partial y^2}+\frac{\partial^2 f_3}{\partial z^2}\right\}\,\mathbf{i}_3.
\end{align}
The previous equation could make us get the wrong idea that the action of the vector Laplacian is the action of the scalar Laplacian component by component, nevertheless this happens only for the case of the Cartesian coordinates and in any other orthogonal curvilinear system the result is totally different as we will see in the next example.

The generalization of the Bitsadze operator to the space $\mathbb{R}^3$ in Cartesian coordinates is:
\begin{align}
\widetilde{\vec{\Delta}}[\vec{f}]&=\left\{\frac{\partial^2 f_1}{\partial x^2}-\frac{\partial^2 f_1}{\partial y^2}-\frac{\partial^2 f_1}{\partial z^2}+2\frac{\partial^2 f_2}{\partial x \partial y}+ 2\frac{\partial^2 f_3}{\partial x \partial z}\right\}\,\mathbf{i}_1\notag\\
&\;\;\;+\left\{-\frac{\partial^2 f_2}{\partial x^2}+\frac{\partial^2 f_2}{\partial y^2}-\frac{\partial^2 f_2}{\partial z^2}+2\frac{\partial^2 f_1}{\partial x \partial y}+ 2\frac{\partial^2 f_3}{\partial y \partial z}\right\}\,\mathbf{i}_2\notag\\
&\;\;\;+\left\{-\frac{\partial^2 f_3}{\partial x^2}-\frac{\partial^2 f_3}{\partial y^2}+\frac{\partial^2 f_3}{\partial z^2}+2\frac{\partial^2 f_1}{\partial x \partial z}+ 2\frac{\partial^2 f_2}{\partial y \partial z}\right\}\,\mathbf{i}_3.
\end{align}

\subsection{Spherical coordinates}
In spherical coordinates:
\begin{align*}
x&=r\sin\theta\cos\psi,\\
y&=r\sin\theta\sin\psi,\\
z&=r\cos\theta,
\end{align*}
where $0\leq r<\infty$, $0\leq\theta\leq\pi$ and $0\leq\psi<2\pi$.

Then the metric coefficients are  $h_1=1$, $h_2=r$ and $h_3=r\sin\theta$; while the unit tangent vectors to the orthogonal coordinate curves are $\mathbf{u}_1=\sin\theta\cos\psi\mathbf{i}_1+\sin\theta\sin\psi\mathbf{i}_2+\cos\theta\mathbf{i}_3$, $\mathbf{u}_2=\cos\theta\cos\psi\mathbf{i}_1+\cos\theta\sin\psi\mathbf{i}_2-\sin\theta\mathbf{i}_3$ and $\mathbf{u}_3=-\sin\psi\mathbf{i}_1+\cos\psi\mathbf{i}_2$.

Hence, the Moisil-Theodoresco operator in spherical coordinates takes the form:
\begin{align}
D_{MT}[f]&=-\left\{\frac{\partial f_1}{\partial r}+\frac{2}{r}f_1+\frac{1}{r}\frac{\partial f_2}{\partial \theta}+\frac{\cot\theta}{r}f_2+\frac{1}{r\sin\theta}\frac{\partial f_3}{\partial \psi}\right\}\notag\\
&\;\;\;\;+\left\{\frac{\partial f_0}{\partial r}+\frac{1}{r}\frac{\partial f_3}{\partial \theta}+\frac{\cot\theta}{r}f_3-\frac{1}{r\sin\theta}\frac{\partial f_2}{\partial \psi}\right\}\,\mathbf{u}_1\notag\\
&\;\;\;\;+\left\{\frac{1}{r}\frac{\partial f_0}{\partial \theta}+\frac{1}{r\sin\theta}\frac{\partial f_1}{\partial \psi}-\frac{\partial f_3}{\partial r}-\frac{1}{r}f_3\right\}\,\mathbf{u}_2\notag\\
&\;\;\;\;+\left\{\frac{1}{r\sin\theta}\frac{\partial f_0}{\partial \psi}+\frac{\partial f_2}{\partial r}+\frac{1}{r}f_2-\frac{1}{r}\frac{\partial f_1}{\partial \theta}\right\}\,\mathbf{u}_3.
\end{align}

The escalar Laplacian is:
\begin{equation}
\Delta_0[f_0]=\frac{\partial^2 f_0}{\partial r^2}+\frac{2}{r}\frac{\partial f_0}{\partial r}+\frac{1}{r^2}\frac{\partial^2 f_0}{\partial \theta^2}+\frac{\cot\theta}{r^2}\frac{\partial f_0}{\partial \theta}+\frac{1}{r^2\sin^2\theta}\frac{\partial^2 f_0}{\partial \psi^2},
\end{equation}
and the vector Laplacian:
\begin{align}
\vec{\Delta}[\vec{f}]&=\left\{\frac{\partial^2 f_1}{\partial r^2}+\frac{2}{r}\frac{\partial f_1}{\partial r}-\frac{2}{r^2}f_1+\frac{1}{r^2}\frac{\partial^2 f_1}{\partial \theta^2}+\frac{\cot\theta}{r^2}\frac{\partial f_1}{\partial \theta}\right.\notag\\
&\;\;\;\;\;+\left.\frac{1}{r^2\sin^2\theta}\frac{\partial^2 f_1}{\partial \psi^2}-\frac{2}{r^2}\frac{\partial f_2}{\partial \theta}-\frac{2\cot\theta}{r^2}f_2-\frac{2}{r^2\sin\theta}\frac{\partial f_3}{\partial \psi}\right\}\,\mathbf{u}_1\notag\\
&\;\;\;+\left\{\frac{\partial^2 f_2}{\partial r^2}+\frac{2}{r}\frac{\partial f_2}{\partial r}-\frac{1}{r^2\sin^2\theta}f_2+\frac{1}{r^2}\frac{\partial^2 f_2}{\partial \theta^2}+\frac{\cot\theta}{r^2}\frac{\partial f_2}{\partial \theta}\right.\notag\\
&\;\;\;\;\;+\left.\frac{1}{r^2\sin^2\theta}\frac{\partial^2 f_2}{\partial \psi^2}+\frac{2}{r^2}\frac{\partial f_1}{\partial \theta}-\frac{2\cot\theta}{r^2\sin\theta}\frac{\partial f_3}{\partial \psi}\right\}\,\mathbf{u}_2\notag\\
&\;\;\;+\left\{\frac{\partial^2 f_3}{\partial r^2}+\frac{2}{r}\frac{\partial f_3}{\partial r}-\frac{1}{r^2\sin^2\theta}f_3+\frac{1}{r^2}\frac{\partial^2 f_3}{\partial \theta^2}+\frac{\cot\theta}{r^2}\frac{\partial f_3}{\partial \theta}\right.\notag\\
&\;\;\;\;\;+\left.\frac{1}{r^2\sin^2\theta}\frac{\partial^2 f_3}{\partial \psi^2}+\frac{2}{r^2\sin\theta}\frac{\partial f_1}{\partial \psi}+\frac{2\cot\theta}{r^2\sin\theta}\frac{\partial f_2}{\partial \psi}\right\}\,\mathbf{u}_3.
\end{align}

We must emphasize here that, unlike what one can think naturally, the components of the previous operator are different than the scalar Laplace operator (there are several extra terms) and this happens for any orthogonal curvilinear coordinate system different from the Cartesian system.

In \cite{Mo} it is shown an example where one can see how the solutions can change if one uses the scalar Laplacian instead of the vector Laplacian for the case of the current distribution for skin effect and for induction heating, one employing the zero-order Bessel function while the other uses the first-order Bessel function.

Finally, generalization of the Bitsadze operator to the space $\mathbb{R}^3$  in spherical coordinates is:

\begin{align}
\widetilde{\vec{\Delta}}[\vec{f}]&=\left\{\frac{\partial^2 f_1}{\partial r^2}+\frac{2}{r}\frac{\partial f_1}{\partial r}-\frac{2}{r^2}f_1-\frac{1}{r^2}\frac{\partial^2 f_1}{\partial \theta^2}-\frac{\cot\theta}{r^2}\frac{\partial f_1}{\partial \theta}\right.\notag\\
&\;\;\;\;\;-\left.\frac{1}{r^2\sin^2\theta}\frac{\partial^2 f_1}{\partial \psi^2}+\frac{2\cot\theta}{r}\frac{\partial f_2}{\partial r}+\frac{2}{r}\frac{\partial^2 f_2}{\partial r\partial \theta}+\frac{2}{r\sin\theta}\frac{\partial^2 f_3}{\partial r\partial\psi}\right\}\,\mathbf{u}_1\notag\\
&\;\;\;+\left\{-\frac{\partial^2 f_2}{\partial r^2}-\frac{2}{r}\frac{\partial f_2}{\partial r}-\frac{1}{r^2\sin^2\theta}f_2+\frac{1}{r^2}\frac{\partial^2 f_2}{\partial \theta^2}+\frac{\cot\theta}{r^2}\frac{\partial f_2}{\partial \theta}\right.\notag\\
&\;\;\;\;\;-\left.\frac{1}{r^2\sin^2\theta}\frac{\partial^2 f_2}{\partial \psi^2}+\frac{2}{r^2}\frac{\partial f_1}{\partial \theta}+\frac{2}{r}\frac{\partial^2 f_1}{\partial r\partial \theta}+\frac{2}{r^2\sin\theta}\frac{\partial^2 f_3}{\partial \theta\partial\psi}\right\}\,\mathbf{u}_2\notag\\
&\;\;\;+\left\{-\frac{\partial^2 f_3}{\partial r^2}-\frac{2}{r}\frac{\partial f_3}{\partial r}+\frac{1}{r^2\sin^2\theta}f_3-\frac{1}{r^2}\frac{\partial^2 f_3}{\partial \theta^2}-\frac{\cot\theta}{r^2}\frac{\partial f_3}{\partial \theta}\right.\notag\\
&\;\;\;\;\;+\left.\frac{1}{r^2\sin^2\theta}\frac{\partial^2 f_3}{\partial \psi^2}+\frac{2}{r^2\sin\theta}\frac{\partial f_1}{\partial \psi}+\frac{2}{r\sin\theta}\frac{\partial^2 f_1}{\partial r\partial \psi}+\frac{2}{r^2\sin\theta}\frac{\partial^2 f_2}{\partial \theta\partial\psi}\right\}\,\mathbf{u}_3.
\end{align}

\section{Applications}

A three-dimensional displacement field $\vec{f}$ in a homogeneous isotropic linear elastic material without volume forces is described by the Lam\'{e}-Navier system:
\begin{equation}\label{Lame}
\mathcal{L}_{\lambda,\mu}[\vec{f}]:=\mu\vec{\Delta}[\vec{f}]+(\mu+\lambda)\mathrm{grad}[\mathrm{div}[\vec{f}]]=0,
\end{equation}
where $\mu>0$, $\displaystyle\lambda>-\frac{2}{3}\mu$. See \cite{Bar}, \cite{La} and \cite{Mi}.

It is worth mentioning here that in \cite{Mi} the operator (\ref{VectorDeltaTilde}) appears when the Poisson constant $\displaystyle \sigma:=\frac{\lambda}{2(\lambda+\mu)}$ equals $\displaystyle \frac{3}{4}$, and for this case the equation:
\[-\widetilde{\vec{\Delta}}[\vec{f}]=0,\]
have an infinity of solutions.

We now proceed obeying the calculation from \cite{MoMoAbBo}. Adding (\ref{VectorDelta}) to (\ref{VectorDeltaTilde}) one gets:
\begin{equation}
\mathrm{grad}[\mathrm{div}[\vec{f}]]=-\frac{1}{2}\left(D_{MT}^2[\vec{f}]+D_{MT}D_{MT}^r[\vec{f}]\right).
\end{equation}

As a matter of fact, the Lam\'{e} equation (\ref{Lame}) can be rewritten in the form:
\begin{equation}
\left(\frac{\mu+\lambda}{2}\right)D_{MT}D_{MT}^r[\vec{f}]+\left(\mu+\frac{\mu+\lambda}{2}\right)D_{MT}^2[\vec{f}]=0,
\end{equation}

Let us denote $\displaystyle\alpha:=\frac{\mu+\lambda}{2}$ and $\displaystyle\beta:=\frac{3\mu+\lambda}{2}$, then we have:
\begin{equation}
\mathcal{L}_{\lambda,\mu}[\vec{f}]=\alpha D_{MT}D_{MT}^r[\vec{f}]+\beta D_{MT}^2[\vec{f}].
\end{equation}
Having in mind the conditions relating $\lambda,\mu$ in (\ref{Lame}), it is easily seen that $\alpha\neq0$ and $\beta\neq0$.

The corresponding expressions for $\mathcal{L}_{\lambda,\mu}[\vec{f}]$ in general orthogonal curvilinear coordinates:
\begin{align}
\mathcal{L}_{\lambda,\mu}[\vec{f}]&=\mu\vec{\Delta}[\vec{f}]+(\mu+\lambda)\mathrm{grad}[\mathrm{div}[\vec{f}]]\notag\\
&=\mu\left(\mathrm{grad}[\mathrm{div}[\vec{f}]]-\mathrm{curl}[\mathrm{curl}[\vec{f}]]\right)+(\mu+\lambda)\mathrm{grad}[\mathrm{div}[\vec{f}]]\notag\\
&=\left\{\frac{2\mu+\lambda}{h_1}\frac{\partial}{\partial q_1}\left[\frac{1}{h_1h_2h_3}\frac{\partial(h_2h_3f_1)}{\partial q_1}\right]+\frac{2\mu+\lambda}{h_1}\frac{\partial}{\partial q_1}\left[\frac{1}{h_1h_2h_3}\frac{\partial(h_1h_3f_2)}{\partial q_2}\right]\right.\notag\\
&\;\;\;\;\;+\frac{2\mu+\lambda}{h_1}\frac{\partial}{\partial q_1}\left[\frac{1}{h_1h_2h_3}\frac{\partial(h_1h_2f_3)}{\partial q_3}\right]-\frac{\mu}{h_2h_3}\frac{\partial}{\partial q_2}\left[\frac{h_3}{h_1h_2}\frac{\partial(h_2f_2)}{\partial q_1}\right]\notag\\
&\;\;\;\;\;+\frac{\mu}{h_2h_3}\frac{\partial}{\partial q_2}\left[\frac{h_3}{h_1h_2}\frac{\partial(h_1f_1)}{\partial q_2}\right]+\frac{\mu}{h_2h_3}\frac{\partial}{\partial q_3}\left[\frac{h_2}{h_1h_3}\frac{\partial(h_1f_1)}{\partial q_3}\right]\notag\\
&\;\;\;\;\;\left.-\frac{\mu}{h_2h_3}\frac{\partial}{\partial q_3}\left[\frac{h_2}{h_1h_3}\frac{\partial(h_3f_3)}{\partial q_1}\right]\right\}\,\mathbf{u}_1\notag
\end{align}
\vspace{-0.5cm}
\begin{align}
\;\;\;\;\;\;\;\;\;\;&\;\;\;+\left\{\frac{2\mu+\lambda}{h_2}\frac{\partial}{\partial q_2}\left[\frac{1}{h_1h_2h_3}\frac{\partial(h_2h_3f_1)}{\partial q_1}\right]+\frac{2\mu+\lambda}{h_2}\frac{\partial}{\partial q_2}\left[\frac{1}{h_1h_2h_3}\frac{\partial(h_1h_3f_2)}{\partial q_2}\right]\right.\notag\\
&\;\;\;\;\;+\frac{2\mu+\lambda}{h_2}\frac{\partial}{\partial q_2}\left[\frac{1}{h_1h_2h_3}\frac{\partial(h_1h_2f_3)}{\partial q_3}\right]-\frac{\mu}{h_1h_3}\frac{\partial}{\partial q_3}\left[\frac{h_1}{h_2h_3}\frac{\partial(h_3f_3)}{\partial q_2}\right]\notag\\
&\;\;\;\;\;+\frac{\mu}{h_1h_3}\frac{\partial}{\partial q_3}\left[\frac{h_1}{h_2h_3}\frac{\partial(h_2f_2)}{\partial q_3}\right]+\frac{\mu}{h_1h_3}\frac{\partial}{\partial q_1}\left[\frac{h_3}{h_1h_2}\frac{\partial(h_2f_2)}{\partial q_1}\right]\notag\\
&\;\;\;\;\;\left.-\frac{\mu}{h_1h_3}\frac{\partial}{\partial q_1}\left[\frac{h_3}{h_1h_2}\frac{\partial(h_1f_1)}{\partial q_2}\right]\right\}\,\mathbf{u}_2\notag
\end{align}
\vspace{-0.5cm}
\begin{align}
\;\;\;\;\;\;\;\;\;\;&\;\;\;+\left\{\frac{2\mu+\lambda}{h_3}\frac{\partial}{\partial q_3}\left[\frac{1}{h_1h_2h_3}\frac{\partial(h_2h_3f_1)}{\partial q_1}\right]+\frac{2\mu+\lambda}{h_3}\frac{\partial}{\partial q_3}\left[\frac{1}{h_1h_2h_3}\frac{\partial(h_1h_3f_2)}{\partial q_2}\right]\right.\notag\\
&\;\;\;\;\;+\frac{2\mu+\lambda}{h_3}\frac{\partial}{\partial q_3}\left[\frac{1}{h_1h_2h_3}\frac{\partial(h_1h_2f_3)}{\partial q_3}\right]-\frac{\mu}{h_1h_2}\frac{\partial}{\partial q_1}\left[\frac{h_2}{h_1h_3}\frac{\partial(h_1f_1)}{\partial q_3}\right]\notag\\
&\;\;\;\;\;+\frac{\mu}{h_1h_2}\frac{\partial}{\partial q_1}\left[\frac{h_2}{h_1h_3}\frac{\partial(h_3f_3)}{\partial q_1}\right]+\frac{\mu}{h_1h_2}\frac{\partial}{\partial q_2}\left[\frac{h_1}{h_2h_3}\frac{\partial(h_3f_3)}{\partial q_2}\right]\notag\\
&\;\;\;\;\;\left.-\frac{\mu}{h_1h_2}\frac{\partial}{\partial q_2}\left[\frac{h_1}{h_2h_3}\frac{\partial(h_2f_2)}{\partial q_3}\right]\right\}\,\mathbf{u}_3.
\end{align}

\section*{Acknowledgements}
The authors were partially supported by Instituto Polit\'ecnico Nacional in the framework of SIP programs and by Fundaci\'{o}n Universidad de las Am\'{e}ricas Puebla, respectively.

\end{document}